\documentstyle[11pt]{article}

 \newtheorem{Theorem}{Theorem}

\newtheorem{Lemma}[Theorem]{Lemma}
\newtheorem{Proposition}[Theorem]{Proposition}

\bibstyle{plain}

\newcommand{\al}{\alpha}
\newcommand{\R}{{\mathbf R}}
\newcommand{\ds}{\displaystyle}
\newcommand{\e}{\varepsilon}

\newcommand{\lra}{\longrightarrow}
\newcommand{\ra}{\rightarrow}
\newcommand{\p}{\partial}
\newcommand{\la}{\lambda}

\newcommand{\N}{{\mathbf N}}

\newcommand{\si}{\sigma}

\newcommand{\ph}{\varphi}

\newcommand{\Ho}{{\mathcal H}}
\newcommand{\Lo}{{\mathcal L}}
\newcommand{\Po}{{\mathcal P}}
\newcommand{\Xo}{{\mathcal X}}
\newcommand{\Yo}{{\mathcal Y}}

\newcommand{\Do}{{\mathcal D}}

\title{\bf Symmetry and monotonicity of least energy solutions}

\date{  }

\author{ Jaeyoung BYEON\footnote{
   Department of Mathematics and PMI, Pohang University of Science and Technology, Pohang, Kyungbuk 790-784,
   Republic of Korea, e-mail: jbyeon@postech.ac.kr},
Louis JEANJEAN\footnote{
   D\'epartement de Math\'ematiques UMR 6623, Universit\'e de Franche-Comt\'e,
   16, Route de Gray,
   25030 Besan\c con, France,
   e-mail: louis.jeanjean@univ-fcomte.fr}
$ \; $  and Mihai MARI\c S\footnote{
   D\'epartement de Math\'ematiques UMR 6623, Universit\'e de Franche-Comt\'e,
   16, Route de Gray,
   25030 Besan\c con, France,
   e-mail: mihai.maris@univ-fcomte.fr}
}

\date{}

\begin{document}

\maketitle

\begin{abstract}
We give a simple proof of the fact that for a large class of quasilinear
elliptic equations and systems the solutions that minimize the
corresponding energy in the set of all solutions are radially
symmetric. We require just continuous nonlinearities and no
cooperative conditions for systems. 
Thus, in particular, our results cannot be obtained by 
using the moving planes method. 
In the case of scalar
equations, we also prove that any least energy solution  has a
constant sign and is
monotone with respect to the radial variable.
Our proofs 
rely on results in \cite{M, BZ} and answer questions from
\cite{brezis-lieb, Li}.

\end{abstract}

\section{Introduction}

We consider the system of partial differential equations
\begin{equation}
\label{1}
 - \mbox{div} (| \nabla u_i|^{p-2} \nabla u_i)
 = g_i (u),
 \qquad i = 1, \dots, m,
\end{equation}
where $ u = ( u_1, \dots, u_m) : \R^N \lra \R^m $, $ 1 < p <
\infty$, $|(y_1, \dots , y_N)|^p = \left( \sum_{j=1}^N y_j ^2
\right)^{\frac p2} $, $g_i (0) = 0 $ and there exists $  G \in
C^1(\R^m \setminus \{ 0 \}, \R) \cap C(\R^m,\R) $ such that
$g_i(u) = \frac{ \p G}{\p u_i} (u) $ for $ u \neq 0$.

\medskip

Formally,
 solutions of (\ref{1}) are critical points of the following energy functional
$$
S(u) = \frac 1p \int_{\R^N} \sum_{i=1}^m |\nabla u_i|^p \, dx -
\int_{\R^N} G(u) \, dx.
$$

The aim of this note is to prove, under general assumptions, that
those solutions of (\ref{1}) which minimize the energy $S$ in the
set of all solutions are radially symmetric (up to a translation
in $\R^N$).  In the scalar case we also study the sign and
monotonicity of these solutions. We do not consider here the
problem of existence of solutions (respectively of least energy
solutions) for (\ref{1}). We believe that our results cover all
situations where the existence of a least energy solution is
already known in the literature.

\medskip

We begin with some definitions. Let $\Pi $ be an affine hyperplane
in $\R^N$,  let $\Pi^+$ and $ \Pi ^-$ be the two closed
half-spaces determined by $\Pi$ and $ s_{\Pi}$ the symmetry
with respect to $ \Pi$
(i.e. $s_{\Pi} (x) = 2 p_{\Pi}(x) - x$, 
where  $p_{\Pi}$ is the orthogonal projection onto $\Pi$).
 Given  a function $f$ defined on $\R^N$,
we define
\begin{equation} \label{2}
f_{\Pi^+}(x) = \left\{
\begin{array}{ll}
f(x) & \mbox{ if }  x \in \Pi^+   \\
f(s_{\Pi} (x)) & \mbox{ if }  x \in \Pi^-
\end{array}
\right. ,
\quad
f_{\Pi^-}(x) = \left\{
\begin{array}{ll}
f(x) & \mbox{ if }  x \in \Pi^- \\
f(s_{\Pi} (x)) & \mbox{ if }  x \in \Pi^+ .
\end{array}
\right.
\end{equation}
For $ \si > 0$, we  denote $f_{\si }(x) = f ( \frac{x}{\si })$. We
say that a  space $ \Xo $ of functions defined on $\R^N$ is
admissible if $\Xo $ is nonempty and
\begin{itemize} \item[{\bf (i)}] $\Xo \subset
L_{loc}^1(\R^N, \R^m) $ and $ \mbox{measure}( \{ x\; | \ \ \;
|u(x)| > \al \}) < \infty $ for any $ u \in \Xo $ and $ \al > 0$;

\item[{\bf (ii)}] $ g_i (u) \in L_{loc}^1(\R^N )$ for any $ u \in
\Xo$ and $ i = 1, \dots,  m$;

\item[{\bf (iii)}] $  \sum_{i=1}^m  |\nabla u_i|^p $ and $ G(u) $
belong to $L^1(\R^N)$ if $ u \in \Xo$;

\item[{\bf (iv)}]$ u_{\si } \in \Xo $ for any $u \in \Xo $ and $
\si
> 0$;

\item[{\bf (v)}]$ u_{\Pi ^+}, \; u_{\Pi ^-}  \in \Xo $ whenever $
u \in \Xo $ and $ \Pi $ is an affine hyperplane in $ \R^N$.
\end{itemize}

\medskip
Let $ \Xo $ be an admissible function space. We note that from
(i) and (iii), $G(0) = 0.$  A function $ u \in \Xo $ is a solution
of (\ref{1}) if it satisfies (\ref{1}) in $\Do '(\R^N)$. If
(\ref{1}) admits solutions in $ \Xo$, we say that $\underline{u} $
is a {\it least energy solution } if $\underline{u} $ is a
nontrivial solution of (\ref{1})  and
\[ S(\underline{u} ) = \inf \{ S(u) \; | \; u \in \Xo \setminus
\{0\}, \  u \mbox{ is a solution of (\ref{1})} \}.\]

\medskip

We introduce the functionals
$$
J(u) = \frac 1p \int_{\R^N} \sum_{i=1}^m  |\nabla u_i|^p \, dx
\qquad
\mbox{ and }
\qquad
V(u) = \int_{\R^N} G(u) \, dx.
$$
Clearly, these functionals are well-defined on any admissible
function space. As we will see, the least energy solutions of
(\ref{1}) come from the following minimization problem:
$$ 
\mbox{ minimize } \; J(u) \;  \mbox{ in the set } \;  \{  \ u \in \Xo \; | \;  V(u ) = \la \}. 
\eqno{(\Po_{\la})}
$$
We shall prove that under some general  conditions 
(see {\bf (C1)-(C3)} or {\bf (D1)-(D3)} below),  
all least energy solutions
of (\ref{1}) in the set $\Xo$ are radially symmetric, up to a
translation in $\R^N$.

\medskip

It is easy to see that $ J(u_{\si }) = \si ^{ N-p} J(u) $ and
$V(u_{\si }) = \si ^N V(u)$. If  $ V(u) > 0$ for some $ u \in
\Xo,$  we have $V(u_{\si }) = 1 $ for $ \si = V(u)^{- \frac 1N}$.
Then, denoting
$$
T = \inf \; \{ J(u) \; | \; u \in \Xo \textrm{ and }  V(u) = 1 \},
$$
we see that
\begin{equation}
\label{3}
 J(v) \geq T \left( V(v) \right)^{\frac{ N-p}{N} }
\qquad \mbox{ for any } v \in \Xo \mbox{ satisfying } V(v) > 0.
\end{equation}
It is  clear that $u$ is a minimizer for problem ($\Po _{\la} $) 
above ($\la > 0$) if and only if $ u _{\si _1}$ is a
minimizer for ($\Po _1$), where $ \si _1 = \la ^{- \frac 1N}$.

\medskip

We assume first that $ 1 < p < N$ and the following conditions are
satisfied.
 \begin{itemize}
\item[\bf{(C1)}]  $T > 0$ and problem ($ \Po _1$) has a minimizer
$ u_* \in \Xo$; \item[\bf{(C2)}] Any minimizer $ u \in \Xo $ of
($\Po _1 $) is a $C^1$ function and satisfies the Euler-Lagrange
system of equations
\begin{equation}\label{4}
 - \mbox{div} (| \nabla u_i|^{p-2} \nabla u_i) = \al g_i (u)
 \qquad \mbox{ in } \Do '(\R^N)
 \end{equation}
for $ i = 1, \dots, m $ and some  $\al \in \R$; \item[\bf{(C3)}]
 Any solution $ u \in \Xo $ of (4) (and not only any minimizer!)
satisfies the Pohozaev identity
\begin{equation}\label{5} (N-p)
J(u) = \al N V(u).  \end{equation}
\end{itemize}

A few comments are in order. Clearly, the most important of the
conditions above is {\bf (C1)}. To our knowledge, the existence of
a minimizer for ($\Po _1 $), under sufficiently general
assumptions on the functions $g_i$  and for arbitrary $ m \in
\N^*$ and $ p \in (1, \infty)$, is still an open problem. However,
several particular cases have been extensively studied in the
literature. A series of papers has been devoted to the case $ p =
2$ and fairly optimal conditions on $ g_i$ that guarantee {\bf
(C1)} have been found by Berestycki-Lions \cite{berestycki-lions}
for $m=1$ and by Brezis-Lieb \cite{brezis-lieb} for $m \geq 1$. In
the case $ m =1$ and $ 1 < p < N$ the existence of a minimizer for
($\Po _1$) has also been proved in \cite{ferrero-gazzola} under
general assumptions on $ g= g_1$ (similar to the assumptions in
\cite{berestycki-lions}). 
Under the conditions  considered in \cite{berestycki-lions} and
\cite{ferrero-gazzola}, the
functionals $J$ and $V$ are well defined on $W^{1,p}(\R^N)$ and
this is clearly an admissible function space.
The setting in \cite{brezis-lieb} also corresponds to our assumptions.

If $T >0$ and ($\Po _1 $) admits minimizers, in most
applications it is quite standard 
to prove that {\bf(C2)} and {\bf (C3)} hold.
 This is indeed the case under the assumptions in
\cite{berestycki-lions, brezis-lieb, ferrero-gazzola}.

\medskip

Next we  consider the case $ p=N $. Note that in this case the
Pohozaev identity (5) becomes $ \al N V(u) = 0$; hence any
``reasonable" solution $u$ of (1) should satisfy $V(u) =0$. 
Since we are interested in  nontrivial solutions, we consider the minimization problem
$$
\mbox{ minimize } \; J(u) \; \mbox{ in the set } \ \;   
 \{  \ u \in \Xo \setminus \{0\} \; | \;  V(u) = 0 \}. \eqno{(\Po _0 ')}
$$

We assume that the following conditions are satisfied.
 \begin{itemize}
\item[{\bf (D1) }] $T_0 := \inf  \{ J(u) \; | \;   u \in \Xo ,  u
\neq 0, V(u) = 0 \} >0$ and ($\Po _0 '$) admits a minimizer $u_0$;
\item[{\bf (D2) }] Any minimizer $ u \in \Xo $ of ($\Po _0 '$) is
$ C^1$ and satisfies the Euler-Lagrange equations (4) for some
$\al
> 0$;

\item[{\bf (D3) }] Any solution $ u \in \Xo $ of (4) (with $\al
> 0$) satisfies the Pohozaev identity $V(u) = 0$.
\end{itemize}

For $p=N=2,$  fairly optimal conditions on $ g_i$ that
guarantee {\bf (D1)-(D3)} have been found by
Berestycki-Gallou\"et-Kavian \cite{bgk} for $m=1$  and by Brezis-Lieb
\cite{brezis-lieb} for $m\ge 1$.

In the next section we show that least energy solutions are
minimizers of $(\Po_{\la })$ for some particular choice of $\la$
if $1 < p < N$, respectively minimizers of $(\Po _0 ')$ if $ p =N$. 
Then we obtain the radial symmetry of such solutions as a
direct consequence of the general results  in  \cite{M}
(in the case  $N=p$, we need some extra-argument in addition to the results in 
 \cite{M}.

In the third section we consider the scalar case $ m =1$ and we
prove that least energy solutions have constant sign and, 
if they tend to zero at infinity, 
then they are monotone with respect
to the radial variable.

In the final section we make some connections with related results
of symmetry and monotonicity in the literature.
 Let us just mention  that,
especially in the scalar case, the symmetry and monotonicity of
solutions of (1) have been studied by many authors, see e.g. 
\cite{gnn, SZ, DPR, BuSi}  and references therein. However, in all these
works it is assumed that the solutions are nonnegative and some
further assumptions on the nonlinearity $g$ are made. They
require, at least, $g$ to be Lipschitz continuous and 
to satisfy a cooperative condition in the case of  systems.
 In the present work, we do not make any
additional assumptions on $g$, except those that guarantee the
existence of least energy solutions (basically, we need $g$ to be
merely continuous and to satisfy some growth conditions near zero,
see \cite{brezis-lieb} and \cite{ferrero-gazzola}). We prove that
our solutions have constant sign and our results are valid as well
for compactly supported solutions and for solutions that do not
vanish. Of course, there is a price we have to pay: our method
works only for least energy solutions, not for any nonnegative
solution of (1).


\section{Variational characterization and symmetry}

We begin with the case $ 1 < p < N$.

\begin{Lemma} \label{l1}  Assume that $ 1 < p < N$ and the
conditions {\bf (C1)-(C3) } hold.
\begin{itemize}
\item[\bf{(i)}] Let $ u$ be a minimizer for ($\Po _1$). Then $
u_{\si _0} $ is a least action solution of (1), where $ \si _0 =
\left( \frac{N-p}{N} T \right)^{\frac 1p}$, and $ S(u_{\si _0}) =
p ( N - p)^{\frac Np - 1} N^{- \frac Np } T ^{ \frac Np }$.
\item[\bf{(ii)}] Let $ v$ be a least energy solution for (1). Then
$v$ is a minimizer for ($\Po _{\la}$), where $ \la = \left(
\frac{N-p}{N} T \right)^{\frac Np}$.
\end{itemize}
 \end{Lemma}
{\bf Proof. } {\bf(i)} By {\bf (C2) }  we know that $ u \in C^1 $
and $u$ satisfies (\ref{4}) for some $ \al \in \R$. Then (\ref{5})
implies $ (N-p) J(u) = \al N V(u)$, which gives $ \al = \frac{
N-p}{N} T > 0$. It is easy to see that $u_{\si _0 } $ satisfies
(1) for $ \si _0 = \al ^{\frac 1p}$ and \[ S( u_{\si _0}) = \si _0
^{ N -p} J(u) - \si _0 ^N V(u) = \si _0 ^{ N -p} T - \si _0 ^N = p
( N - p)^{\frac Np - 1} N^{- \frac Np } T ^{ \frac Np }.\]
 Let $ w\in \Xo$, $ w \neq 0$, be a solution of (\ref{1}). By {\bf (C3) }
we have $ ( N-p) J(w) = N V(w)$. If $ J(w) = 0$, we have $ \nabla
w = 0 $ a.e. on $ \R^N$, hence $w$ must be constant. Since $\mbox{
measure}\{ x \in \R^N \; | \; | w(x) > \al \} < \infty $ for any $
\al > 0$, we infer that $ w = 0 $, a contradiction.
Thus $ J(w)
> 0$ and $ V(w) = \frac{ N-p}{N} J(w) > 0$. On the other hand, by
(\ref{3}) we get $J(w) \geq T \left( V(w)\right)^{\frac{N-p}{N}}
$, i.e. $ J(w) \geq T \left(  \frac{ N-p}{N} J(w)
\right)^{\frac{N-p}{N}} $, which gives
\begin{equation}\label{6}
J(w) \geq \left( \frac{ N-p}{N} \right)^{\frac{ N-p}{p}} T^{\frac
Np}.
\end{equation}
Combined with Pohozaev identity, this implies
\begin{equation}\label{7}
 S(w) = J(w) - V(w)
 = \frac pN J(w) \geq
p ( N - p)^{\frac Np - 1} N^{- \frac Np } T ^{ \frac Np } = S(u_{\si_0})
\end{equation}
and we infer that $u_{\si_0} $ is a least energy solution for
(\ref{1}).

{\bf(ii)} Conversely, let $v$ be a least energy solution for
(\ref{1}). Then $ (N-p) J(v) = N V(v)$ by {\bf (C3)}, hence $ S(v)
= \frac pN J(v)$. It is obvious that the inequalities (\ref{6})
and (\ref{7}) above  are satisfied with $w = v.$ On the other
hand, $ S(v) = S(u_{\si _0})$ and we infer that $v$ must satisfy
(7) with equality sign, that is, \[ J(v) = \left( \frac{ N-p}{N}
\right)^{\frac{ N-p}{p}} T^{\frac Np} \ \textrm { and  } \  V(v) =
\frac {N-p}{N} J(v) = \left( \frac{N-p}{N}\right)^{\frac Np}
T^{\frac Np}.\] A simple scaling argument shows that $v$ is a
minimizer for ($\Po _{\la }$), where $ \la = \left(
\frac{N-p}{N}\right)^{\frac Np} T^{\frac Np}$; equivalently,
$v_{\si _1} $ is a minimizer for ($\Po _1$), where $ \si _1 =
\left( \frac{ N-p}{N} T \right)^{ - \frac 1p} = \si _0 ^{-1}$.
This completes the proof of Lemma 1. \hfill $\Box $

The symmetry of least energy solutions will follow  from Lemma
\ref{l1} and a general symmetry result in \cite{M}. For the
convenience of the reader, we recall here that result.

 \begin{Theorem}[\cite{M}]   Assume that $ u : \R^N \lra
\R^m $ belongs to some function space $\Yo $ and solves
 the minimization problem
$$
\begin{array}{l}
\mbox{ minimize } \ds  \int_{\R^N} F( u(x),  |\nabla u(x)|) \, dx \ 
\\
 \mbox{ in the set } \; 
\Big \{  u \in \Yo \; \Big | \; 
\ds   \int_{\R^N} H( u(x), |\nabla u(x)|) \, dx = \la\neq 0 \Big \}.
\end{array}
\eqno{(\Po )}
$$
Suppose that the following conditions are satisfied:
\begin{itemize}
\item[\bf{(A1)}] For any $ v \in \Yo $ and any affine hyperplane
$\Pi$ in $\R^N$
 we have $ v_{\Pi ^+}, v_{\Pi ^-} \in \Yo$.
\item[\bf{(A2)}] Problem ($\Po$) admits minimizers in $\Yo$ and
any  minimizer is a $C^1$ function on $\R^N$.
\end{itemize}
Then, after a translation, $u$ is radially symmetric.
\end{Theorem}

Lemma 1 implies that least energy solutions solve the minimization
problem ($\Po _{\la }$) for some $ \la >0$. Conditions {\bf (C1)},
{\bf (C2) } and   property {\bf (v)}  in the definition of admissible
spaces imply that ($\Po _{\la }$)  satisfies the assumptions of
Theorem 2. Thus we get:

 \begin{Proposition} \label{p1}
Assume that $1<p<N$ and {\bf (C1)-(C3) } hold. Then (1) admits a
least energy solution and each least energy solution is radially
symmetric (up to a translation in $\R^N$).
 \end{Proposition}

Now we turn our attention to  the case $ p=N $.

 \begin{Proposition} \label{p2} Assume that $ p = N $ and {\bf (D1)-(D3) } hold. Then
(1) admits a least energy solution and any least energy solution
solves ($\Po _0 '$).

Moreover, if we assume that $G$ is either negative or positive in
some ball $B_{\R^m} (0, \e ) \setminus \{ 0\}$ and $u\in \Xo $ is
a least energy solution such that $ u(x) \lra 0 $ as  $|x| \lra
\infty $, then $u$  is radially symmetric (up to a translation in
$\R^N$).
 \end{Proposition}

\medskip

\noindent
{\bf Proof. } Let $u_0 $ be a minimizer for ($\Po _0 '$). By {\bf
(D2)} and {\bf (D3)} we have  $ V(u_0) = 0 $ and $ u_0$ satisfies
(4) for some $ \al >0$. Let $ u_1 = (u_0)_{\si }$, where $ \si =
\al ^{\frac 1p}$. It is easy to see that $ u_1 $ solves (1) and $
S(u_1) = J(u_1) - V(u_1) = J(u_0) - \si ^N V(u_0) = J(u_0) = T_0$.
For any solution $ u \in \Xo$, $ u \neq 0$ of (1) we have $ V(u) =
0 $ by {\bf (D3) } and $ S(u) = J(u) \geq T_0 = J(u_1) $. Hence $
u_1 $ is a least energy solution.

If $ v$ is a least energy solution, then $ V(v) = 0$  by {\bf (D3)
} and $J(v) = S(v) = S(u_1) = T_0$, thus $v$ solves ($\Po _0 '$).

Although Theorem 2 does not apply directly to minimizers of
problem ($\Po _0 '$) (because the value of the constraint in ($\Po
_0 '$) is zero), its proof can still be adapted to those
minimizers. Indeed, the proof of Theorem 2 shows that whenever $u$
is a minimizer of ($\Po $) and $ \Pi $ is an affine hyperplane
such that $u_{\Pi ^+}$ and $u_{\Pi ^-}$ are also minimizers, $u$
must be symmetric with respect to $\Pi$. The only place where the
assumption $ \la \neq 0$ is used in Theorem 2 is to show that for
any $e \in S^{N-1}$ there exists an affine hyperplane $ \Pi $
orthogonal  to $e$
such that 
\begin{equation}\label{8}
\int_{\Pi ^-} H(u(x), |\nabla u(x)|) \, dx = \int_{\Pi ^+} H(u(x),
|\nabla u(x)|) \, dx =\frac{ \la}{2}.
\end{equation}
From (\ref{8}) it follows then easily that $u_{\Pi ^+}$ and
$u_{\Pi ^-}$ are also minimizers.

In the present case we will use the fact that $G(u)$ has a
constant sign in a neighborhood  of $\infty$ to find hyperplanes
that ``split the constraint in  two equal parts." A similar idea has
already been used in \cite{LM}. Henceforth we assume that $u$ is a least action
solution, $ u(x) \lra  0 $ as $ | x| \lra \infty $ and, say,
$G(\xi ) < 0$ for $0 < |\xi | < \e $.
For $ e \in S^{N-1}$ and $ t \in \R$, we denote $\Pi_{e, t} = \{ x
\in \R^N \; | \; x \cdot e = t \}$, $\Pi_{e, t}^- = \{ x \in \R^N
\; | \; x \cdot e < t \}$ and $\Pi_{e, t}^+ = \{ x \in \R^N \; |
\; x\cdot e > t \}$. We claim that for any $ e \in S^{N-1}$, there
exists $ t_e \in \R$ such that
\begin{equation}\label{9}
\ds \int_{\Pi_{e, t_e}^-} G(u(x))\, dx = \ds \int_{\Pi_{e, t_e}^+}
G(u(x))\, dx = 0 \qquad \mbox{ and } \qquad u_{\Pi_{e, t_e}^-}
\not\equiv 0, \; u_{\Pi_{e, t_e}^+} \not\equiv 0.
\end{equation}
To see this, fix $ e \in S^{N-1}$ and define $ \ph _e^\pm(t) = \ds
\int_{\Pi_{e, t}^\pm} G(u(x))\, dx,$ respectively. It follows
that $ \ph_e^+$  and $ \ph_e^-$ are continuous because $ G(u) \in L^1(\R^N)$.
Since $u$ is continuous, $ u \not\equiv 0$,  $\lim_{|x| \to \infty}u(x) = 0$ 
and $ G < 0$ on $ B_{\R^m} (0 , \e) \setminus \{ 0 \}$,  
it is not hard to see that 
there exist $t^-, t^+ \in \R$, $ t^- < t^+$ such that
$$
\ph_{e}^-(t^-) < 0, \quad \ph_{e}^+(t^+) < 0 \quad \textrm { and  } \quad
u_{\Pi^-_{e,t^-}} \neq 0, \quad u_{\Pi^+_{e,t^+}} \neq 0.
$$
 Since
$\ph_{e}^+(t^-) = V(u) - \ph_{e}^-(t^-)= - \ph_{e}^-(t^-)$, it follows  that
 $\ph_{e}^+(t^+) < 0 < \ph_{e}^+(t^-).$   
From the mean value property, we see that there exists $t_e \in (t^-, t^+)$
 satisfying (\ref{9}). It is clear  that
$ u_{\Pi_{e, t_e}^-} , \; u_{\Pi_{e, t_e}^+} \in \Xo \setminus \{0 \}$ 
because $ \Xo $ is admissible
and (9) implies that $V(u_{\Pi_{e, t_e}^-}) = V(u_{\Pi_{e, t_e}^+}) = 0$, hence
$ J( u_{\Pi_{e, t_e}^-}) \geq T_0$, 
$ J( u_{\Pi_{e, t_e}^+}) \geq T_0$.
 On the other hand, it is easy to see that 
$J( u_{\Pi_{e, t_e}^-})+  J( u_{\Pi_{e, t_e}^+}) = 2 J(u) = 2 T_0$. 
Thus $ J( u_{\Pi_{e, t_e}^-}) =   J( u_{\Pi_{e, t_e}^+}) = T_0$ and 
$u_{\Pi_{e, t_e}^-} $, $ u_{\Pi_{e, t_e}^+}$ are also minimizers for ($\Po _0 '$).
Then arguing exactly as in the proof of Theorem 2 in \cite{M}, it follows that 
after a translation, $u$ is radially symmetric.  
 \hfill $\Box$

\section{Monotonicity results}

Throughout this section we assume that $ m =1$. We consider the
following additional conditions for an admissible space $\Xo.$
 \begin{itemize}
\item[{\bf (vi)}] 
For any $ u \in \Xo $ and $ t \geq 0$, $ s \leq 0$, we
have $\min(u, t) \in \Xo$ and $ \max(u, s) \in \Xo$.  
\item[{\bf (vii)}]
If $ u \in \Xo$ and $ u \geq 0$ (respectively $u \leq 0$), then $
u^* \in \Xo$ (respectively $- (-u)^* \in \Xo$), where $ u^*$ is
the Schwarz rearrangement of $u$.
\end{itemize}

\begin{Proposition} \label{signn}
Let  an admissible space $\Xo$ satisfy the  condition
{\bf(vi)}.  Assume that $1<p<N$ and {\bf (C1)} holds. 
If $u \in \Xo$ is a solution of ($\Po _{\la }$) for some $ \la >0$, then
$u$ does not change sign.
\end{Proposition}

\noindent{\bf Proof.} This is a simple consequence of scaling. Indeed, let $
u_+ = \max(u, 0)$ and $ u_- = \min (u, 0)$. It is clear that $ V(u_+ ) + V(u_-)
= V(u) = \la $ and $ J(u_+ ) + J(u_-) = J(u)$. If $ V(u_-) < 0$, then
necessarily $ V(u_+) > \la$. For $ \si = \left( \frac{\la  }{V(u_+) }
\right)^{\frac1N} \in (0, 1) $ we have $V((u_+)_{\si } ) = \si ^N V(u_+) = \la
$ and $J((u_+)_{\si } ) = \si ^{N-p} J(u_+) \leq \si ^{N-p} J(u) < J(u)$,
contradicting the fact that $u$ is a minimizer. Thus necessarily $ V(u_-) \geq
0$. In the same way $ V(u_+) \geq 0$, therefore $ V(u_-), V(u_+) \in [0, \la]$.
Using inequality (3) (which trivially holds if $V(v) = 0$), we get
$$
T\la ^{\frac{N-p}{N}} = J(u) =  J(u_+ ) + J(u_-) \geq TV(u_+)^{\frac{N-p}{N}} +
TV(u_-)^{\frac{N-p}{N}},
$$
which gives
\begin{equation} \label{10}
1 \geq \left(\frac{V(u_+)}{\la }\right)^{\frac{N-p}{N}} +
\left(\frac{V(u_-)}{\la }\right)^{\frac{N-p}{N}}.
\end{equation}
Since $ V(u_+ ) + V(u_-) = \la$, (\ref{10}) implies that either $
V(u_+) = 0 $ or $ V(u_-) = 0 $. If $ V(u_-) = 0 $ and $ V(u_+) =
\la $ we see that $u_+ $ satisfies the constraint and
\begin{equation} \label{11}
J(u_+) = J(u) - J(u_-) \leq J(u).
\end{equation}
Since $u$ is a minimizer, we must have equality in (11) and this
gives $J(u_-) = 0$, hence $ u_- = 0$ and $ u = u_+ \geq 0$.
Similarly $ V(u_+) = 0$ implies $ u = u_- \leq 0$. \hfill $\Box$

\begin{Proposition} \label{critical}
Let  an admissible space $\Xo$ satisfy the condition  {\bf(vi)}.  
Assume that $p=N$ and {\bf (D1)} holds. We have:
\begin{itemize}
\item[{ \bf (a)}] 
if $ G < 0$ on $[-\e, 0) \cup(0, \e]$ for some
$ \e >0$,  then $ u \in \Xo $ is a minimizer of ($\Po _0 '$) if and only
if it solves the problem
$$
\mbox{ minimize } \; J(v)  \; \mbox{ in the set } 
 \{  \  v \in \Xo \; | \; v \neq 0, V(v)\geq 0  \}; \eqno{(\Po _0 '')}
$$
 \item[{\bf (b)}] if $ G > 0$ on $[-\e, 0) \cup(0, \e]$, then $ u \in \Xo $ solves ($\Po _0
'$) if and only if it solves the problem
$$
\mbox{ minimize } \; J(v)  \; \mbox{ in the set } 
 \{  \  v \in \Xo \; | \; v \neq 0, V(v) \leq 0  \} . \eqno{(\Po _0 ''')}
$$
\end{itemize}
 Moreover, any
minimizer of ($\Po _0 ''$) or ($\Po _0 '''$) does not change sign.
\end{Proposition}

\noindent{\bf Proof.} It clearly suffices to prove {\bf (a)}.

Consider $ v \in \Xo $ such that $ v \geq 0$ a.e. and $ V(v) > 0$. For
$ t\geq 0$ we define $v^t(x) = \min(v(x), t)$. By {\bf(vi)} we
have $ v^t\in \Xo$. We claim that there exists $ t_* > 0$ such
that $ V(v^{t_*} ) = 0$.

The continuity of $G$, properties {\bf (i)} and {\bf (iii)} in the definition
of admissible spaces and the dominated convergence theorem imply
that the mapping $ t \mapsto V(v^t) = \ds \! \int_{\R^N} \! G(v^t(x))
\, dx $ is continuous on $(0, \infty)$. 
Since $G(v^{ \e } (x)) < 0$ whenever $ v(x) \neq 0$ and we cannot have 
$ v (x) = 0 $ a.e. because $V(v) > 0$, we infer that  $V(v^{ \e }) < 0$.

We claim that there exists $ t_0 > \e $ such that $ V(v^{t_0}) >
0$. Two situations may occur:

\smallskip

Case 1.  There exists an increasing sequence $ t_n \to \infty $ 
such that $\{G(t_n)\}_{n=1}^\infty$ is bounded from below. 
Let $ m =  \inf_{n \geq 1} G(t_n)$. 
By dominated convergence we get
$$
V(v^{t_n}) - V(v) = \ds \int_{\{v \geq t_n\} } G(t_n) - G(v(x)) \,
dx \geq \int_{\{ v \geq t_n\} } m - G(v(x) ) \, dx  \lra 0 \quad
\mbox{ as } n \lra \infty;
$$
hence $ V(v^{t_n}) \geq \frac 12 V(v) > 0$ for $n$ sufficiently large.

\smallskip

Case 2.  $G(s) \lra - \infty $ as $  s \lra \infty$. Then,
since $ v \geq 0 $ a.e. and $ V(v) > 0$, we see that the set $ A =
\{ s >0 \; | \; G(s) > 0 \}$ is nonempty. Let $M = \sup A < \infty$. 
It follows that $G(s ) \leq 0$ for $ s \ge M$. It is
clear that $M > \e$ and $ V(v^M) \geq V(v) > 0$.
The claim is thus proved. 

\smallskip

Now the continuity of the mapping $ t \longmapsto V(v^t)$ 
implies that there exists 
$ t_* \in (\e , t_0)$ such that $V(v^{t_*}) = 0$. 
Similarly, if $w\in \Xo$, $w\leq 0$ a.e. and $ V(w) > 0$ there is some 
$\tilde{t} >0$ such that $V(- (-w)^{\tilde{t}}) = 0$.

Next let $ u_0 \in \Xo $ be a minimizer of ($\Po _0'$). Suppose $
V(u) > 0$ for some $ u \in \Xo. $
Then at least one of the
quantities $V(u_+)$ and $V(u_-)$ is positive. If $V(u_+) >0$, take
$ t_* > 0$ such that $V(u_+^{t_*}) = 0$. We have $ u_+^{t_*} \in
\Xo \setminus \{ 0 \}$ and
\begin{equation} \label{12}
J(u) \geq J(u_+) \geq J( u_+^{t_*} ) \geq J(u_0) = T_0.
\end{equation}
Hence $\ds \inf \{ J(u) \; | \; u \in \Xo, u \neq 0, V(u) \geq 0
\} = J(u_0) = T_0 $ and $u_0$ is a solution of  ($\Po _0''$).

Conversely, assume that $ u $ is a solution of ($\Po _0''$). We prove that
\begin{equation} \label{13}
V(u_+) = V(u_-) = V(u) = 0.
\end{equation}
We argue again by contradiction. If (\ref{13}) does not hold, the
inequality  $ V(u_+) + V(u_-) = V(u) \geq 0$ implies that at least
one of the quantities $V(u_+)  $ and $V(u_-)$ must be positive.
Suppose that $V(u_+)
>0$. As above we find $ t_* > 0$ such that $V(u_+^{t_*}) = 0$ and
then (12) holds for $u$. Moreover, since $u$ is a minimizer of
($\Po _0''$) we have $ J(u) \le T_0$ and therefore all
inequalities in (12) are in fact equalities. But $J(u_+) = J(
u_+^{t_*} ) $ implies $\ds \int_{\{ u > t_* \} } |\nabla u|^p \,
dx = 0$, hence $ \nabla u = 0 $ a.e. on $\{ u > t_* \} $ which
gives $ \nabla (( u - t_*)_+) =0$ a.e. and we infer that $( u -
t_*)_+ = 0$ a.e., that is $ u \leq t_* $ a.e. Then we have $u_+ =
u_+^{t_*}$ and consequently $V(u_+) = V( u_+^{t_*}) = 0$, contrary
to our assumption. We argue similarly if  $V(u_-)  >0$ and
(\ref{13}) is proved. Since $ V(u) = 0 $ and $ J(u) = T_0 =
J(u_0)$, we see that $u$ solves ($\Po _0'$).

Lastly we show that if $u$ is a minimizer of ($\Po _0''$),
then  either $u_+ = 0 $ a.e. or $ u_- = 0 $ a.e. 
(but we cannot have $ u_+  = u_- = 0$ a.e. because $J(u) = T_0 > 0$). 
Indeed, if $u^+ \neq 0$ and $ u^- \neq 0$, (\ref{13}) would imply 
$ J(u_+) \geq T_0 $ and  $ J(u_-) \geq T_0 $ and this would give
 \[T_0=J(u) =  J(u_+) + J(u_-) \geq 2T_0 > 0,\] which is a contradiction.
This completes the proof.
\hfill $ \Box$
 
\medskip

Next we prove  the monotonicity of scalar minimizers.

\begin{Theorem} \label{sign}
Let $ \Xo$ be an admissible space satisfying the conditions
{\bf(vi)} and {\bf(vii)}.
We assume that  conditions {\bf (C1)-(C3)} hold if $ 1 < p < N$, 
respectively conditions {\bf (D1)-(D3)} hold if $ p = N$.
In the case  $p = N$, we also assume that 
 there exists $ \e > 0$ such that either $ G>0$ or $ G < 0$ on
 $[-\e, 0) \cup (0, \e]$.
Then any least energy solution $u$ of (\ref{1}) such that  
$\lim_{|x| \to \infty}u(x) = 0$ is, up to a translation, radially symmetric
and  monotone with respect to $r = |x| \in [0, \infty)$.
\end{Theorem}

\noindent{\bf Proof.} 
Symmetry follows directly from Propositions \ref{p1} and \ref{p2}. 
Hence there is a function $ \tilde{u} : [0, \infty ) \lra \R$ such that 
$ u(x) = \tilde{u} ( |x|) = \tilde{u} (r)$.  
From Lemma \ref{l1} and Proposition \ref{p2}, we know that any 
least energy solution is a  minimizer of ($\Po _{\la }$) for some $ \la > 0$,  
respectively of ($\Po _0'$). 
We will show that whenever $ u(x) = \tilde{u} (r) $ 
solves one of these minimization problems and tends to zero at infinity, 
$\tilde{u}$ is monotone on $[0, \infty)$. 

We have to introduce some notation. In what follows, $ \Ho ^{N-1}$
is the $(N-1)-$dimensional Hausdorff measure and $\Lo ^N$ is the
Lebesgue measure  on $\R^N$. Given a Lebesgue measurable set
$E\subset \R^N$, we denote by $\p ^* E$ its measure theoretic
boundary, i.e. $\p ^* E = \{ x \in \R^N \; | \; 0 < {\ds \lim_{r \ra
0}} \frac{\Lo ^N ( E \cap B(x, r) )}{\Lo ^N( B(x, r))} < 1 \}.$ If
$E$ has finite measure, we denote by $E^*$ the Schwarz
rearrangement of $E$, i.e. $E^*$ is the open ball centered at the
origin such that $\Lo^N(E) = \Lo ^N(E^*)$. We recall the
isoperimetric inequality: if $ E \subset \R^N$ is bounded and
measurable, then 
\begin{equation}\label{14} 
\Ho^{N-1} ( \p (E^*)) \leq \Ho ^{N-1}( \p ^* E),
\end{equation} 
with equality holding if and
only if $E$ is equivalent to a ball in $\R^N$ (see, e.g.,
Proposition 2.2 p. 157 in \cite{BZ} and references therein). Note
that the right side in (\ref{14}) might be $\infty$.

Now let  $u$ be as above. From
Proposition \ref{signn}   and Proposition \ref{critical}, we know 
that $u$ has constant sign;  hence we may
assume that $u \ge 0$. Let  $u^*$
be the Schwarz rearrangement of $u$. For $ t \geq 0$ we denote
$$
E_t = \{ x \in \R^N \; | \; u(x) > t \}, \qquad F_t = \{ x \in
\R^N \; | \; u^*(x) > t \}.
$$ 
Note that $E_t$ is bounded for any $ t > 0$ because $u$ tends
to zero at $\infty$ and  $(E_t)^* = F_t$.

We argue by contradiction and we assume that  $ \tilde{u} $ is not
nonincreasing. Then there exist $ 0 \leq r_1 < r_2$ such that 
$ 0< \tilde{u}(r_1) < \tilde{u}(r_2)$. Since $ \tilde{u}(x) \lra 0 $ 
as $ |x| \lra \infty$, there
exists $ r_3 > r_2$ such that $u(r_3) =u(r_1)$. Denoting 
$ a = u(r_1)$ and $ b = u(r_2)$, we see  that for any $ t \in ( a, b)$,
$E_t$ is nonempty and is not equivalent to a ball.
The isoperimetric inequality gives 
\begin{equation} 
\label{15}
\Ho ^{N-1} ( \p F_t) \leq \Ho ^{ N-1} ( \p ^* E_t) \qquad \mbox{ for
any } t \in (0, M), 
\end{equation} 
with {\it strict } inequality for $ t \in (a, b)$.

Since $u \in C^1$ and $\lim_{|x| \to \infty}u|x) = 0,$ we see 
that $u$ is bounded and 
$  (u  -t)_+ ,  (u ^* -t)_+ \in W^{1, p}(\R^N)$ for any $ t > 0$. 
Let $M = \max_{x \in \R^N}u(x)$.
Using the
coarea formula for $W^{1, p}$ functions (see, e.g., Proposition
2.1 p. 157 in \cite{BZ}), we find
\begin{equation} \label{16}
\int_{\{ u \geq t \} } |\nabla u|^p \, dx = \int_t^M \left(
\int_{u^{-1}(s)} |\nabla u |^{p-1} \, d \Ho ^{N-1} \right) \, ds.
\end{equation}
The coarea formula for $  (u ^* -t)_+$ gives
\begin{equation} \label{17}
\int_{\{ u ^* \geq  t \} } |\nabla u^* |^p \, dx = \int_t^M
\left(| \nabla u^*( (u^*)^{-1}(s)) |^{p-1} \Ho ^{N-1}
((u^*)^{-1}(s)) \right) \, ds.
\end{equation}
Passing to the limit as $ t \downarrow 0$ and using the monotone
convergence theorem, we see that (\ref{16}) and (\ref{17}) also
hold for $t =0$.

The following result is a simple consequence of Lemma 3.1 p. 161 in \cite{BZ}.

\begin{Proposition} \label{brothers} {\bf (\cite{BZ})}
Let $ v \in \Xo $ be a nonnegative function 
that tends to zero at infinity and let $ v^* (x) = \tilde{v}^*  (|x|) $ 
be the Schwarz rearrangement of $v$. 
There exists a set $ N_v \subset (0, \sup (v))$ of Lebesgue measure zero such that 
for any $ t \in (0, \sup (v)) \setminus N_v\, $, 
$ (\tilde{v}^*)^{-1}(t) $ contains only  one point,
$(\tilde{v}^*)' ((\tilde{v}^*)^{-1}(t)) $ exists, $\Ho ^{N-1} ( v^{-1}(t)) $ and
$\Ho ^{N-1} ( (v^*)^{-1}(t)) $ are finite and 
\begin{equation} \label{18}
\int_{v^{-1}(t) } |\nabla v|^{p-1} \, d\Ho ^{N-1} \geq |\nabla
v^*((v^*)^{-1}(t)) |^{p-1} \Ho ^{N-1} ((v^*)^{-1}(t)).
\end{equation}
Moreover, if $ t \in (0, \sup (v)) \setminus N_v $ and we have
equality in (18) then necessarily
\begin{equation} \label{19}
\Ho^{N-1} ( \p^* \{ x \in \R^N \; | \; v(x ) > t \} ) = \Ho ^{
N-1} ( v^{-1}(t)) =  \Ho ^{ N-1} ( (v^*)^{-1}(t))
\end{equation}
and $ | \nabla v| = |\nabla v^* ( (v^*)^{-1}(t)) | = constant \; $
$\Ho^{N-1}- $a.e. on $ v^{-1}(t)$.
\end{Proposition}

By Proposition \ref{brothers} we infer that $u$ satisfies
(\ref{18}) for any $ t \in (0, M) \setminus N_u$, where $ \Lo ^1(
N_u) = 0$. Moreover, the isoperimetric inequality (\ref{15})
(which is strict for $ t \in (a, b))$ implies that $u$ cannot
satisfy (\ref{19}) for $ t \in (a, b)\setminus N_u$. Therefore we
have {\it strict } inequality in (\ref{18}) for $u$ whenever $ t
\in (a,b) \setminus N_u$. Integrating (18) from $0$ to $M$ and
using (\ref{16}) and (\ref{17}) (with $ t =0$) we get
 \begin{equation} \label{20}
 \ds \int_{\R^N} |\nabla u |^p \, dx > \int_{\R^N} |\nabla u ^*|^p \, dx,
\mbox{ or equivalently } J(u) > J(u^*).
 \end{equation}
On the other hand it is clear that  $ u^* \in \Xo \setminus \{
0\}$
 and   $ V(u^*) = V(u)$, therefore (20) contradicts
the fact that $u$ is a minimizer.
This  proves that $\tilde{u}$ must be nonincreasing.
\hfill $ \Box $

\section{Some  remarks and examples}

\noindent
{\bf Remark 9 } 
In the scalar case $m=1$ it is well known (see for
example the Introduction of \cite{Br}) that if $g$ is odd then any
least energy solution has a constant sign. In Remark II.6 of
\cite{Li}, Lions raised the question (for $p=2$ and $N\geq 3$) whether
this remains true without assuming $g$ odd. Proposition \ref{signn}
gives an affirmative answer for any $1<p <N$ and Proposition
\ref{critical}, under some mild additional assumptions, for $p=N$.
Previous partial results were obtained by Brock \cite{Br}, using
rearrangement arguments, assuming that $ 1 < p \leq 2$,  
the minimizer $u$ satisfies $u(x) \lra 0$ as $|x| \to \infty$ and 
$g \in C^{0,p-1}(\R)$. Nothing was proved for $p>2$. 

\medskip

\noindent
{\bf Remark 10 } 
If $ N \geq 3$, $ p =2$, $ m =1$ and under the 
assumption that $g$ is odd, the existence of least energy solutions
for (1) has been proved in \cite{berestycki-lions} by showing that
problem $(\Po _1)$ admits a minimizer. The minimizer found in
\cite{berestycki-lions} was radial by construction, but it was not
known whether all least energy solutions were radially symmetric.
The existence of a minimizer for $(\Po _1)$ without the oddness 
assumption on $g$ has also been proved in  \cite{Li},  but nothing 
was known about the symmetry or the  sign of such minimizers.
Our results imply that any least energy solution is radially symmetric, 
has constant sign and is monotone with respect to  the radial variable, 
no matter whether $g$ is odd or not.

\smallskip

In the case $ N \geq 2$, $ p =2$, $ m \in \N^*$, the existence of least 
energy solutions  is also known  (see \cite{brezis-lieb} for general
results, historical notes, comments and further references). 
If $N > 2$, the existence of a minimizer for ($\Po _{\la }$) and  the
existence of least energy solutions have been proved in
\cite{brezis-lieb} under very general assumptions on the functions
$g_i$. It has also been shown that the solutions are smooth
(Theorem 2.3 p. 105 in \cite{brezis-lieb}) and satisfy the
Pohozaev identity (Lemma 2.4 p. 104 in \cite{brezis-lieb}).
However, as already mentioned in \cite{brezis-lieb} p. 99, the
existence of radially symmetric least energy solutions was not
clear. Indeed, the Schwarz symmetrization that lead to a radial
minimizer in \cite{berestycki-lions} could not be used in
\cite{brezis-lieb} because of the general assumptions on the
nonlinearity made there. In fact, it is known that the Schwarz
rearrangements may be used for systems only if the nonlinearity
satisfies a cooperative condition.

Proposition \ref{p1}  above implies that all least energy solutions
 of the system considered  in \cite{brezis-lieb}
 are radially symmetric.

\smallskip

If $ N=2 $ and $ G(\xi ) < 0 $ for $ 0 < |\xi | \leq \e $, the
existence of least energy solutions and the existence of
minimizers for ($\Po _0'$) have  been proved in
\cite{bgk,brezis-lieb}. It has also been shown that such solutions
are smooth, satisfy the Pohozaev identity and tend to 0 as $ | x |
\lra \infty$. Therefore Proposition \ref{p2} implies that any
least energy solution is radially symmetric.

\smallskip

We have to mention that if  $p=2$ and if the minimizers of ($\Po
_{\la }$) satisfy a unique continuation principle, it has already
been proved in \cite{lop1} that any minimizer is radially
symmetric (modulo translation). In \cite{lop1} no cooperative
condition is required when $m \geq 2$ but using a unique
continuation principle require in particular $g$ to be $C^1$. Our
results are still valid when a unique continuation principle fails
(e.g., for minimizers with compact support). Note that compactly
supported minimizers may occur in some applications (cf. Theorem
3.2 (ii) p. 111 in \cite{brezis-lieb}; see also \cite{M} for such
an example). In the scalar case $m=1$, \cite{lop1} does
not say anything about the sign of the minimizers.


\medskip

\noindent
{\bf Remark 11 } 
If $ 1< p < N$ and $ m =1$, it has been
proved in \cite{ferrero-gazzola}, under general conditions
on $g$,  that  problem ($\Po _{\la }$) 
admits minimizers (thus (1) has least energy solutions). The
minimizers found in \cite{ferrero-gazzola} were radially symmetric
by construction. It follows from Proposition \ref{p1}  that any
least energy solution is radially symmetric.

If,  in addition to the assumptions of Theorem \ref{sign}, it is 
assumed that  $g$ is locally Lipschitz on $(0, \infty)$ and non-increasing
on some interval $ (0, s_0)$ and $1<p<2$,  it has been proved in
\cite{DPR}  that any nonnegative  solution of (1) is
radially symmetric and that $ u(x) = \tilde{u}(|x|) $ satisfies 
$\tilde{u} ' (r) < 0$ whenever $ r>0$ and $ \tilde{u}(r) > 0$. 
The same result is true when $p>2$ if it is assumed in addition 
 that the critical set of the solution $u$ is
reduced to one point (see  \cite{SZ}).
These assumptions are not necessary for us but,  of course,  we only
deal with least energy solutions. 


\medskip

\noindent
{\bf Remark 12 } 
(i) The symmetry results in Section 2 hold without
any change if we replace the functional $J$ by a functional of the
form $ \ds \int_{\R^N} \sum_{i =1}^m A_i(u, \nabla u_i) dx$ where
$\xi \to A_i(u, \xi)$ is $p$-homogeneous for any $i=1,...,m$.

\medskip

(ii) Our method still works for more general functionals of the
form
$$
\tilde{J}(u) = \frac 1p \int_{\R^N} |x|^{\al } \sum_{i=1}^m A_i(u) |\nabla u_i|^p \, dx
\qquad
\mbox{ and }
\qquad
\tilde{V}(u) = \int_{\R^N} |x|^{\beta} G(u) \, dx.
$$
In this case, using Theorem 1 in \cite{M},
we obtain that minimizers (and the corresponding minimum action solutions) are
axially symmetric.

Functionals of this type appear, e.g., in the
Caffarelli-Kohn-Nirenberg problem (which consists in minimizing
$\ds \int_{\R^N} |\nabla u|^q |x|^{- aq} \; dx $ under the
constraint $ \ds \int_{\R^N} |u|^p |x|^{-bp } \, dx = const.$,
where $ q>1$, $ p>1$, $ a \leq b < \frac Nq$ and $0< \frac 1q -
\frac 1p = \frac{1+a-b}{N}$). It has  been proved that minimizers
for this problem exist and, in general, are not radially symmetric
(see \cite{bw} and references therein).


\bigskip

{\bf Acknowledgment} The research of the first author was
 supported  in part by KRF-2007-412-J02301 of Korea Research Foundation.


\begin{thebibliography}{}

\bibitem{berestycki-lions} {\sc H. Berestycki, P.-L. Lions, }
{\it Nonlinear scalar field equations, I. Existence of a ground
state, } Arch. Rational Mech. Anal. 82 (1983),  313-345.

\bibitem{bgk}
{\sc H. Berestycki, T. Gallou\"et,  O. Kavian,} {\it Equations
de champs scalaires euclidens non lin\'{e}aires dans le plan,} 
C.R. Acad. Sc. Paris S\'erie I - Math.  297  (1983), 307-310
and Publications du Laboratoire  d'Analyse Num\'erique,
Universit\'{e} de Paris VI, 1984.

\bibitem{brezis-lieb} {\sc H. Br\'ezis, E. H. Lieb, }
{\it Minimum Action Solutions for Some Vector Field Equations, }
Comm. Math. Phys. 96 (1984), 97-113.

\bibitem{bw}{\sc J. Byeon,  Z.-Q. Wang,}
{\it Symmetry breaking of extremal functions for the
Caffarelli-Kohn-Nirenberg inequalities,} Comm. Contemp. Math. 4
(2002), 457--465.

\bibitem{Br}{\sc F. Brock,}
{\it Positivity and radial symmetry of solutions to some
variational problems in $\R^N$, }   J. Math. Anal. Appl. 296 (2004),
226--243.

\bibitem{BZ} {\sc J. E. Brothers, W.P. Ziemer, }
{\it Minimal rearrangements of Sobolev functions, } J. reine
angew. Math. 384 (1988),  153-179.

\bibitem{BuSi} {\sc J. Busca, B. Sirakov, }
{\it Symmetry results for semilinear elliptic systems in the whole
space, } J. Diff. Eq. 163, No. 1 (2000), 41-56.

\bibitem{DPR} {\sc L. Damascelli, F. Pacella, M. Ramaswamy, }
{\it Symmetry of ground states of p-Laplace equations via the
moving plane method, } Arch. Rational Mech. Anal. 148 (1999),
291-308.

\bibitem{ferrero-gazzola} {\sc A. Ferrero, F. Gazzola, }
{\it On subcriticality assumptions for the existence of ground
states of quasilinear elliptic equations, } Adv. Diff. Eq. 8, No.
9 (2003),  1081-1106.

\bibitem{gnn}
{\sc  B. Gidas, W. N. Ni,  L. Nirenberg,} 
{\it Symmetry of positive solutions of nonlinear elliptic equations in $\R^n$, }
Adv. Math. Supp. Stud. 7A (1981), 369-403.


\bibitem{Li} {\sc P.L. Lions, }
{\it  The concentration-compactness principle in the Calculus of
Variations, The locally compact case, Part 2, } Ann. Inst. H.
Poincar\'{e}  Anal. Non Lin\'{e}aire 1 (1984), 223-283.


\bibitem{lop1} {\sc O. Lopes, }
{\it  Radial symmetry of minimizers for some translation and
rotation invariant functionals, } J. Diff. Eq. 124 (1996),
378-388.

\bibitem{LM} {\sc O. Lopes,  M. Montenegro, }
{\it Symmetry of mountain pass solutions for some vector field
equations, } J. Dyn. Diff. Eq. 18, No. 4 (2006),  991-999.

\bibitem{M} {\sc M. Mari\c s, }
{\it On the symmetry of minimizers, }
Arch. Rational Mech. Anal., to appear;
 arXiv:0712.3386 (www.arxiv.org)


\bibitem{SZ} {\sc J. Serrin, H. Zou, }
{\it Symmetry of ground states of quasilinear elliptic equations,
} Arch. Rational Mech. Anal. 148 (1999), 265-290.




\end{thebibliography}
\end{document}